\documentclass[abstracton,a4paper]{scrartcl}

\usepackage[headsepline]{scrpage2}
\usepackage[latin1]{inputenc}
\usepackage{amsfonts}
\usepackage{amsmath}
\usepackage{amssymb}
\usepackage{amsthm}

\usepackage{wasysym}  

\usepackage{url}

\usepackage{graphicx}

\usepackage[pdftex,colorlinks,breaklinks,linkcolor=black,citecolor=black,filecolor=black,menucolor=black,urlcolor=black,pdfauthor={Andreas Weber},pdftitle={Analysis of the physical Laplacian on a locally finite graph}, plainpages=false, pdfpagelabels, bookmarksnumbered=true]{hyperref}

\ihead{A. Weber: Analysis of the  physical Laplacian and the heat flow on a locally finite graph}

\newtheorem{theorem}{Theorem}[section]
\newtheorem{lemma}[theorem]{Lemma}
\newtheorem{proposition}[theorem]{Proposition}
\newtheorem{corollary}[theorem]{Corollary}

\theoremstyle{definition}
\newtheorem{definition}[theorem]{Definition}

\theoremstyle{remark}
\newtheorem{remark}[theorem]{\bf Remark}

\DeclareMathAlphabet{\lier}{U}{eur}{m}{n}  

\newcommand{\N}{\mathbb{N}}

\newcommand{\R}{\mathbb{R}}

\newcommand{\C}{\mathbb{C}}

\newcommand*\e{\mathrm{e}}

\newcommand*\supp{\mathrm{supp}}

\newcommand*\dom{\mathrm{dom}}

\newcommand*\pddt{\frac{\partial}{\partial t}}

\newcommand*\D{{\cal D}}

\title{Analysis of the physical Laplacian and the heat flow on a locally finite graph}
\author{ Andreas Weber\footnote{
						   E-mail: andreasweber.mail@gmail.com,
						   Address: Institut f\"ur Algebra und Geometrie,
						   Universit\"at Karlsruhe (TH),
						   Englerstr. 2, 76128 Karlsruhe, Germany.}\\
						    Universit\"at Karlsruhe (TH)}

\date{}

\begin{document}

\maketitle 

\begin{abstract}
	We study the physical Laplacian and the corresponding heat flow on an infinite, locally finite graph
	with possibly unbounded valence.\\
	
         						
\end{abstract}		

\section{Introduction and preliminaries}

Our aim in this paper is to study the physical Laplacian and the corresponding heat flow on a
connected, locally finite graph $G=(V,E)$. In contrast to the the normalized Laplacian 
(sometimes also called  combinatorial Laplacian), the physical Laplacian is not always a bounded 
operator on $\ell^2(V)$ and hence, its analysis is more complicated.
While the normalized Laplacian has been studied extensively in the past 
(cf.  \cite{MR1421568,0645.05048} and the references therein) investigations concerning the
physically more motivated unbounded (physical) Laplacian have started just recently, cf.
\cite{MR2257129,Keller:2008nr,MR2199204,wojciechowski-2007,Wojciechowski:yq2008}, and as pointed out in \cite{Keller:2008nr} the spectral properties of these Laplacians might
be very different. 
Note also that spectral properties of the physical Laplacian on locally tessellating planar graphs 
are studied in \cite{Keller:rt,Keller:yq}. 

We first show in the next two section that the Laplacian $\Delta$ with appropriate domain is a positive, essentially self-adjoint operator on $\ell^2(V)$, cf. Theorem \ref{essential self-adjointness} and Proposition \ref{positive}. This is an analogue of the well known result that the Laplacian initially defined
on the set of smooth functions with compact support on a complete Riemannian manifold $M$ extends
to an unbounded self-adjoint operator on $L^2(M)$, cf. \cite{MR0062490}.

In Section \ref{heat equation} we study the heat equation $\pddt u +\Delta u=0$ on $G$.
Similar to the case of the heat equation on non-compact Riemannian
manifolds (cf. \cite{MR711862}), we construct a fundamental solution for infinite graphs by
using an exhaustion of the graph by a sequence of finite subsets of vertices. We also address the question of uniqueness of bounded solutions of the heat equation with respect to some initial condition
$u_0:V\to\R$. In fact, we give in Theorem \ref{unique solution} a condition, which can be interpreted as a weak curvature bound, that ensures the uniqueness of bounded solutions, and we give an example for a graph with unbounded valence with such a weak curvature bound. This generalizes a result by 
J. Dodziuk who proved uniqueness in the case of bounded valence 
(see \cite{MR2246774,MR2218014}).\\

After this work was finished we learned about the recent  work of Radoslaw K. Wojciechowski which contains related results (see \cite{wojciechowski-2007,Wojciechowski:yq2008}). Furthermore, the essentail self-adjointness of the Laplacian was independently proved by Palle Jorgensen in \cite{Jorgensen:2008sf}. See also the preprint \cite{Jorgensen:yg} by Palle Jorgensen and Erin Pearse.\\
In the preprint \cite{Keller:rm}, Matthias Keller and Daniel Lenz extended some of the mentioned results to the more general context of regular Dirichlet forms on discrete sets.\\
For related topics as random walks and analysis on networks, we refer to \cite{Lyons:eu,MR1324344,MR1743100} and the references therein.\\

From now on, we always consider a non-oriented, countable, locally finite, connected graph $G=(V,E)$
with counting measure. 
Furthermore, we denote by $m(x) = \#\{y\in V : y\sim x\}$ the valence of $x\in V$ and by
$$ \ell^2(V) = \{ f: V\to \C ~|~ \sum_{x\in V} |f(x)|^2 < \infty \}$$
with inner product
$$ \langle f,g\rangle = \sum_{x\in V} f(x)\overline{g(x)}$$
the complex  Hilbert space of square summable functions.
Sometimes we also  will need the set of oriented edges 
$E_0 = \{ [x,y], [y,x] : x,y \in V,\, x\sim y \}$. Basically, every edge $e\in E$ is represented by 
two oriented edges in $E_0$. We will also use the notation $[x,y] = -[y,x]$.

The (physical) Laplacian $\Delta$ is  the linear operator defined by
$$ \Delta f (x) = m(x)f(x) - \sum_{y\sim x}f(y) = \sum_{y\sim x}\Big(f(x)-f(y)\Big).$$
Note, that the normalized
Laplacian $\Delta_N$, defined by $\Delta_N f(x) = \frac{1}{m(x)}\Delta f(x)$, is easily seen
to be a bounded operator on the Hilbert space 
$$\ell^2(V,m) = \{  f: V\to \C ~|~ \sum_{x\in V} m(x)|f(x)|^2 < \infty \}.$$

\section{Essential self-adjointness}

We denote by $C_c(V)\subset \ell^2(V)$ the dense subset of functions $f:V\to \C$ with finite support.
Furthermore, we will need the subset
 $$\D = \{ f\in \ell^2(V) : \Delta f \in \ell^2(V) \}$$
which is dense in $\ell^2(V)$ since it contains  $C_c(V)$.

To make the proofs in this section more readable, we define
\begin{eqnarray*}
 \Delta_m:  &C_c(V) \to \ell^2(V),& f\mapsto \Delta f\\
 \Delta_M:  &\D \to \ell^2(V),& f\mapsto \Delta f.
\end{eqnarray*}

\begin{theorem}\label{essential self-adjointness} 
    The operator
	$$ \Delta_m: C_c(V) \to \ell^2(V)$$
	is essentially self-adjoint.			
\end{theorem}

To prove this theorem, we will need the following lemma.
\begin{lemma}
  $\Delta_m: C_c(V) \to \ell^2(V)$ is symmetric.
\end{lemma}
 \begin{proof}
  We have to show for $f,g \in C_c(V)$ that $\langle \Delta_m f,g\rangle = \langle f, \Delta_m g\rangle$.
  Because of
   $$
   	\langle \Delta_m f,g\rangle = \sum_{x\in V} m(x)f(x)\overline{g(x)} -
   						\sum_{x\in V}\sum_{y\sim x}f(y)\overline{g(x)}
   $$						
  and
   $$	
   	\langle f, \Delta_m g\rangle = 	 \sum_{x\in V} m(x)f(x)\overline{g(x)} -
						\sum_{x\in V}\sum_{y\sim x} f(x)\overline{g(y)}
   $$
   this is equivalent to					
   $$
   	\sum_{x\in V}\sum_{y\sim x}f(y)\overline{g(x)} =	 \sum_{x\in V}\sum_{y\sim x} f(x)\overline{g(y)},
   $$					
   which always holds true.
 \end{proof}

The proof of Theorem \ref{essential self-adjointness} will follow from the next proposition and lemmas.

\begin{proposition}[cf. \cite{MR0493419}] 
The symmetric operator $\Delta_m: C_c(V) \to \ell^2(V)$ is essentially self-adjoint
 if and only if
 $\ker(\Delta_m^* \pm i)=\{0\}$. 
\end{proposition}

\begin{lemma}
 The adjoint operator $\Delta_m^*$ of $\Delta_m: C_c(V)\to \ell^2(V)$ is
 $$\Delta_M: \D\to \ell^2(V).$$
\end{lemma}
\begin{proof}
 For any $g\in C_c(V), f\in \ell^2(V)$ we have
 $$ \langle \Delta_m g, f\rangle = \langle g, \Delta f\rangle$$
 and hence, if $\Delta f \in \ell^2(V)$, we obtain $f\in \dom(\Delta_m^*)$, i.e.
 $\D\subset  \dom(\Delta_m^*)$.\\
 Let on the other hand $f\in  \dom(\Delta_m^*)$. Then there is an $h\in \ell^2(V)$ such that
 for any $g\in C_c(V)$
 $$
    \langle \Delta_m g, f\rangle = \langle g, h\rangle.
 $$
 As the left hand side coincides with $\langle g, \Delta f\rangle$ and the set $C_c(V)$  is dense
 in $\ell^2(V)$, we obtain 
 $$
   \Delta f = h = \Delta_m^* f
 $$
 and therefore $\dom(\Delta_m^*)\subset \D$ and $\Delta_m^* = \Delta_M$.
\end{proof}

 \begin{lemma}\textup{(Maximum principle for subharmonic functions, cf. \cite[Lemma 1.6]{MR743744}).}
 \label{max princ for subharmonic}
 Let $f: V\to \R$ satisfy
 $$ \Delta f \leq 0$$
 and assume that there is an $x\in V$ with $f(x) = \max \{ f(y) : y\in V\}$. 
 Then $f$ is constant. 
 \end{lemma}
 \begin{proof} 
  From $ \Delta f(x) \leq 0$   it follows immediately
  $$ m(x) f(x) \leq \sum_{y\sim x} f(y).$$
  Since $f$ attains its maximum at $x$, it follows $f(y)= f(x)$ for any $y\sim x$. As we assume our 
  graph to be connected, the result follows by induction.
 \end{proof}

 \begin{lemma}
  We have 
  $$ \ker( \Delta_M \pm i) = \{0\}.$$
 \end{lemma}
 \begin{proof}
  Let $f\in \ell^2(V)$ such that $\Delta_Mf = if$. Then it follows $(m(x) -i) f(x) = \sum_{y\sim x} f(y)$ and
  therefore
  $$ (m^2(x) +1)^{1/2}\, |f(x)| \leq \sum_{y\sim x} |f(y)|.$$
  This yields
  \begin{eqnarray*}
   \Delta |f|(x)  &=& m(x) |f(x)| - \sum_{y\sim x} |f(y)|\\
   		     &\leq&  (m^2(x) +1)^{1/2}\,  |f(x)| - \sum_{y\sim x} |f(y)|\\
		     &\leq &0.
  \end{eqnarray*}
  Since we assume $f\in\ell^2(V)$, the function $|f|$ attains its maximum and from the 
  maximum principle for subharmonic functions it follows $|f|= const.$ and hence, $f=0$.
  The same proof works for $(\Delta_M + i)$.
 \end{proof}

Hence, the operator  $\Delta_m: C_c(V)\to \ell^2(V)$ is essentially self-adjoint and has therefore a unique self-adjoint extension which we denote in the following by $\bar{\Delta}: \D\to \ell^2(V)$.

\begin{remark}
 The essential self-adjointness of the Laplacian contrasts the fact that the 
 {\em adjacency matrix} $A: C_c(V)\to \ell^2(V),$
 $$ Af(x) = \sum_{y\sim x} f(y)$$
 is in general not essentially self-adjoint if the graph has unbounded valence. A first example for this fact was given by M\"uller in \cite{MR898555}.
 Furthermore,  for any $n\in \N$ there is an infinite graph with deficiency index $n$, cf. 
 \cite[Section 3]{0645.05048} and the references given therein. In the very recent preprint \cite{Golenia:lq} by Gol\'enia this topic is discussed further.
 In the case of bounded valence however, $A: \ell^2(V)\to \ell^2(V)$ is always a bounded self-adjoint operator
 as $A= M - \Delta$ with $Mf(x)= m(x)f(x)$ and both $M$ and $\Delta$ are bounded self-adjoint operators.
\end{remark}

\begin{proposition} 
 Let $G=(V,E)$ denote a locally finite, connected graph. Then the Laplacian $\bar{\Delta}$ is 
 a  bounded operator on $\ell^2(V)$ if and only if the valence is bounded:
 $$ \sup_{x\in V} m(x) < \infty.$$
\end{proposition}
\begin{proof}
 If the valence $m$ is bounded from above, a straightforward calculation using the triangle inequality and the Cauchy-Schwarz inequality leads to
 $$ ||\bar{\Delta}|| \leq 2\sup_{x\in V} m(x).$$
 On the other hand, if $m$ is unbounded we choose a sequence $(x_j)_{j\in \N}$ in $V$ with\linebreak
 $\sup_{j\in\N}m(x_j)=\infty$ and define $f_j: V\to \C$ by $f_j(x_j)=1$ and $f_j(x)=0$ if $x\neq x_j$.
 Then we clearly have $f_j\in\D$ and 
 \[ \Delta f_j(x) = \left\{
    \begin{array}{ll}
      m(x_j), & x=x_j\\
      -1,         & x\sim x_j\\
      0,		& \mbox{else}.
    \end{array}\right.
 \] 
 Hence, $||\Delta f_j||^2 = m(x_j)^2 + m(x_j)$ is unbounded but $||f_j||=1$.
\end{proof}

\section{Co-boundary operator and positivity}

For the set of oriented edges $E_0$  we define
$$ \ell^2(E_0) = \{ \phi: E_0\to \C ~|~ \phi(-e) = -\phi(e),\, \sum_{e\in E} |\phi(e)|^2 < \infty \}.$$
Together with the inner product
$$ (\phi,\psi) = \frac{1}{2}\sum_{e\in E_0} \phi(e)\overline{\psi(e)}$$
$\ell^2(E_0)$ is a Hilbert space.

\begin{definition}[\cite{MR1421568}]
 The map
 $$ d: C_c(V) \to \ell^2(E_0),\, f\mapsto df$$
 with 
 $$ df([x,y]) = f(x)-f(y)$$ 
 is called  {\em co-boundary operator} of the graph $G=(V,E)$.
\end{definition}

\begin{proposition}\textup{(cf. also \cite[Lemma 1.8]{MR743744}).}\label{positive}
  For all $f,g \in C_c(V)$ we have
  $$ (df, dg) = \langle \Delta f, g\rangle$$
  and hence, $\bar{\Delta}$ is positive.
\end{proposition}
 \begin{proof}
  For any oriented edge $e\in E_0$ we denote by $i(e)$, resp. $t(e)$, the initial, resp. terminal, vertex
  of $e$. Then we have for $f, g\in C_c(V)$:
  \begin{eqnarray*}
   (df, dg) &=& \frac{1}{2}\sum_{e\in E_0}\Big(f(i(e)-f(t(e)) \Big)\Big(\overline{g(i(e))}-\overline{g(t(e)}\Big)\\
   	       &=& \frac{1}{2}\sum_{x\in V}\sum_{y\sim x} \Big(f(x)- f(y)\Big)
	       			\Big(\overline{g(x)}-\overline{g(y)} \Big)\\	
	       &=& \frac{1}{2} \sum_{x\in V} (\Delta f)(x)\overline{g(x)} - 
	       		\frac{1}{2}\sum_{x\in V}\sum_{y\sim x} \Big(f(x)- f(y)\Big)\overline{g(y)}.
  \end{eqnarray*}
  A straightforward calculation now shows
  $$ \sum_{x\in V}\sum_{y\sim x} \Big(f(x)- f(y)\Big)\overline{g(y)} = 
  		-\sum_{x\in V}(\Delta f)(x)\overline{g(x)} $$
  and the result follows.		
 \end{proof}

\section{Heat equation}\label{heat equation}

In this section we study the heat equation 
$$ \pddt u + \Delta u=0$$
on the graph $G=(V,E)$. We say that
a function $p: (0,\infty)\times V\times V\to \R$ is a fundamental solution of the heat equation,
if for any bounded initial condition $u_0:V\to \R$, the function 
$$ u(t,x) = \sum_{y\in V} p(t,x,y)u_0(y),\qquad t>0, x\in V$$
is differentiable in $t$, satisfies the heat equation, and if for any $x\in V$
$$\lim_{t\to 0^+}u(t,x) = u_0(x)$$
holds.\\
In Section \ref{heat kernel on an infinite graph} below we construct on any locally finite graph 
a fundamental solution by using an idea similar to the one in the setting of Riemannian manifolds, cf. \cite{MR711862}.  Such a construction was independently developed in
Radoslaw  Wojciechowski's PhD thesis \cite{wojciechowski-2007}.

\subsection{Maximum principles}

For any subset $U\subset V$ we denote by $\mathring{U}=\{ x\in U : y\sim x \Rightarrow y\in U\}$
the interior  of $U$. The boundary  of $U$ is $\partial U=U\setminus \mathring{U}$.
\begin{theorem}\label{maximum principle 1}
 Let $U\subset V$ be finite and $T>0$. Furthermore, we assume that the function
 $u: [0,T]\times U\to \R$ is differentiable with respect to the first component and satisfies on
 $[0,T]\times \mathring{U}$ the inequality
 $$ \pddt u + \Delta u \leq 0.$$
 Then the function $u$ attains its maximum on the parabolic boundary
 $$ \partial_P\Big([0,T]\times U \Big)= \Big(\{0\}\times U\Big)\cup\Big([0,T]\times \partial U \Big).$$
\end{theorem}
 \begin{proof}
  In a first step we assume that $u$ satisfies the strict inequality
  $$ \pddt u + \Delta u < 0$$
  and that at  the point $(t_0,x_0)\in (0,T]\times \mathring{U}$ the function $u$ attains its maximum. 
  Then it follows
  $\pddt u(t_0,x_0)\geq 0$ and hence 
  \begin{eqnarray*}
     0  &>& \Delta u(t_0,x_0) \\
     	&=& \sum_{y\sim x_0} \Big(u(t_0,x_0)-u(t_0,y)\Big).
  \end{eqnarray*} 
  This contradicts $u(t_0,x_0)\geq u(t_0,y)$ for $y\sim x_0$.\\
  
  In the general case, we consider for any $\varepsilon >0$ the function
  $$ v_{\varepsilon}(t,x) = u(t,x) - \varepsilon t.$$
  Then we have
  $$
   \pddt v_{\varepsilon} + \Delta v_{\varepsilon} = \pddt u + \Delta u - \varepsilon < 0.
  $$ 
  Using our first step, we obtain
  \begin{eqnarray*}
   \max \{ u(t,x) : t\in [0,T], x\in U\} &\leq & \max\{v_{\varepsilon}(t,x) : t\in [0,T], x\in U \} + \varepsilon T\\
   		& = & \max\{v_{\varepsilon}(t,x) : (t,x)\in \partial_P([0,T]\times U)  \} + \varepsilon T\\
		& \leq &\max\{ u(t,x) : (t,x)\in \partial_P([0,T]\times U)  \} + \varepsilon T\\
		& \to & \max\{ u(t,x) : (t,x)\in \partial_P([0,T]\times U)  \} \qquad (\varepsilon \to 0).
  \end{eqnarray*}
 \end{proof}
 
 If we assume $U$ to be connected, we can say more:
 \begin{proposition}
   Let $U\subset V$ be finite and connected and $T>0$. Furthermore, we assume that the function
  $u: [0,T]\times U\to \R$ is differentiable with respect to the first component and satisfies on
   $[0,T]\times \mathring{U}$ the inequality
  $$ \pddt u + \Delta u \leq 0.$$
  If $u$ attains its maximum at $(t_0,x_0)\in (0,T]\times \mathring{U}$ we have
  $$ u(t_0,\cdot) = u(t_0,x_0).$$
 \end{proposition}
 \begin{proof}
   Assume that at  the point $(t_0,x_0)\in (0,T]\times \mathring{U}$ the function $u$ attains its maximum.
  Then it follows
  $\pddt u(t_0,x_0)\geq 0$ and hence 
  \begin{eqnarray*}
     0  &\geq& \Delta u(t_0,x_0) \\
     	&=& \sum_{y\sim x_0} \Big(u(t_0,x_0)-u(t_0,y)\Big).
  \end{eqnarray*} 
  But as the difference $u(t_0,x_0)-u(t_0,y)$ is always non-negative we may conclude that
  $u(t_0,y) = u(t_0,x_0)$ for any $y\sim x_0$ and since $U$ is connected, the claim follows.
\end{proof}

A special case of the preceding proposition is the following corollary.
\begin{corollary}\label{maximum principle 3}
  Let $U\subset V$ be finite and connected and $u: U\to \R$ satisfies on $\mathring{U}$ the inequality
  $$ \Delta u \leq 0.$$
  If  $u$ attains its maximum in $\mathring{U}$, the function $u$ is constant.
\end{corollary}

\subsection{Heat equation on domains}

In this subsection $U\subset V$ denotes always a finite subset. We consider the Dirichlet problem (DP)
$$\left\{\begin{array}{lclr}
  \pddt u (t,x) + \Delta_U u(t,x) &=& 0,          & x\in\mathring{U}, t>0\\
  u(0,x)				 &= &u_0(x), & x\in \mathring{U}\\
  u|_{[0,\infty)\times \partial U} &=& 0 &
\end{array}\right.$$
on $U$, where $\Delta_U$ denotes the Dirichlet Laplacian on $\mathring{U}$, i.e.
$$ \Delta_U = \pi\circ \Delta\circ \iota$$
where $\iota: \ell^2(\mathring{U})\to \ell^2(V)$ denotes the canonical embedding and
$\pi: \ell^2(V)\to \ell^2(\mathring{U})$ denotes the orthogonal projection of $\ell^2(V)$ onto
the {\em subspace} $\ell^2(\mathring{U})\subset \ell^2(V)$.

As $\Delta_U: \ell^2(\mathring{U})\to \ell^2(\mathring{U})$ is positive (Proposition \ref{positive}), 
self-adjoint (Theorem \ref{essential self-adjointness}), and 
$\dim \ell^2(\mathring{U}) < \infty$, there are finitely many eigenvalues
$$ 0\leq \lambda_0 \leq \lambda_1 \leq \dots\leq \lambda_n$$
with a corresponding orthonormal basis consisting of real eigenfunctions
$$ \Phi_0, \Phi_1,\ldots, \Phi_n.$$

\begin{lemma} 
 Let $U\subset V$ be finite and $\Delta_U$ the Dirichlet Laplacian on $U$.
 Then there are no non-trivial harmonic functions on  $U$, in particular $\lambda_0>0$.
\end{lemma}
\begin{proof}
 This follows immediately from the maximum principle in Corollary \ref{maximum principle 3}.
\end{proof}

\begin{lemma} \label{heat kernel finite}
The {\em heat kernel} $p_U$ of $\Delta_U$ with Dirichlet boundary conditions
is given by
$$ p_U(t,x,y) = \sum_{j=0}^n \e^{-\lambda_j t}\Phi_j(x)\Phi_j(y).$$
\end{lemma}

\begin{proof}
 This follows immediately from the facts 
  $p_U(t,x,y) = \e^{-t\Delta_U}\delta_y(x),$
 $\e^{-t\Delta_U}\Phi_j = \e^{-t\lambda_j}\Phi_j,$ and
 $\delta_y(x) = \sum_{j=1}^n \langle \Phi_j, \delta_y\rangle \Phi_j$.
\end{proof}

\begin{theorem} For $t>0, x,y\in \mathring{U}$ we have
  \begin{itemize}
   \item[\textup{(a)}] $p_U(t,x,y)\geq 0,$
   \item[\textup{(b)}] $ \sum_{y\in \mathring{U}} p_U(t,x,y) \leq 1$,
   \item[\textup{(c)}] $\lim_{t\to 0^+}\sum_{y\in\mathring{U}}p_U(t,x,y) = 1$,
   \item[\textup{(d)}] $\pddt p_U(t,x,y) = -\Delta_{(U,y)} p_U(t,x,y)$.
  \end{itemize}
\end{theorem}
\begin{proof}
(a) and (b) are immediate consequences of the maximum principle 
(cf. Theorem \ref{maximum principle 1}) and (d) follows from Lemma  \ref{heat kernel finite}. 
For the proof of (c) we remark that this follows from the continuity of the semigroup $\e^{-t\Delta_U}$
at $t=0$ if the limit is understood in the $\ell^2$ sense. As $U$ is finite all norms are equivalent and pointwise convergence follows also.
\end{proof}

\subsection{Heat kernel on an infinite graph}\label{heat kernel on an infinite graph}

Let $U_k \subset V, k\in \N$ be a sequence of {\em finite} subsets with
$U_k\subset \mathring{U}_{k+1}$ and $\bigcup_{k\in \N}U_k = V$. Such a sequence always
exists and can be constructed as a sequence $U_k = B_k(x_0)$ of metric balls with center 
$x_0\in V$ and radius $k$. The connectedness of our graph $G$ implies that the union 
of these $U_k$ equals $V$.

In the following, we will write
$p_k$ for the heat kernel $p_{U_k}$ on $U_k$, and consider $p_k(t,x,y)$ as a function
on $(0,\infty)\times V\times V$ by defining it to be zero if either $x$ or $y$ is not contained in
$\mathring{U_k}$.
Then, the maximum principle implies the monotonicity of the heat kernels, i.e.  
$$   p_k\leq p_{k+1},$$
and the following limit exists (but could be infinite so far).
\begin{definition}  For any $t>0, x,y\in V$, we define
			    $$ p(t,x,y) = \lim_{k\to\infty} p_k(t,x,y).$$	
\end{definition}

From the properties of $p_k$ we immediately obtain
\begin{lemma} For any $t>0, x,y\in V$ we have:
  \begin{itemize}
    \item[\textup{(a)}] $p(t,x,y) = p(t,y,x),$
    \item[\textup{(b)}] $p(t,x,y)\geq 0$.
  \end{itemize}
\end{lemma}

Our aim is to show that $p$ is a fundamental solution (the heat kernel) of the heat equation
on our graph $G=(V,E)$.
For this, we first prove the following proposition.

\begin{proposition}\label{sequence of solutions}
 Let $u_k: (0,\infty)\times V\to \R, k\in \N,$ be a non-decreasing sequence with 
 $\supp\, u_k(t,\cdot)\subset \mathring{U}_k$ such that
 \begin{itemize}
  \item[\textup{(i)}] $\pddt u_k(t,x) = -\Delta_{U_k}u_k(t,x),$
  \item[\textup{(ii)}] $ |u_k(t,x)| \leq C <\infty,$
  			 for some constant $C>0$ that neither depends on $x\in V, t>0$ nor on $k\in \N$.
 \end{itemize}
 Then the limit 
 $$ u(t,x) = \lim_{k\to\infty}u_k(t,x)$$
 is finite and $u$ is a solution for the heat equation. Furthermore, the convergence is uniform on
 compact subsets of $(0,\infty)$.
\end{proposition}

\begin{proof}
 The finiteness of $u(t,x)$ follows from the second assumption.\\
 From Dini's theorem we may conclude that for any $x\in V$ the sequence $u_k(\cdot,x)$
 converges uniformly on compact subsets of $(0,\infty)$ and therefore, the limit $u(\cdot,x)$
 is continuous.\\  
 Furthermore, we have
 \begin{eqnarray*}
   u_k'(t,x)  &=& -\Delta_{U_k} u_k(t,x)\\
   		&=& \left\{\begin{array}{ll} 
			-m(x)u_k(t,x) + \sum_{y\sim x} u_k(t,y), & \mbox{if~} x\in\mathring{U}_k\\
			0, & \mbox{else}
			\end{array}\right.\\
		&\to& -m(x)u(t,x) + \sum_{y\sim x} u(t,y)  = -\Delta u(t,x),
 \end{eqnarray*}
 where the convergence is uniform on compact subsets of $(0,\infty)$.\\
 Hence, the limit $u(\cdot,x)$ is differentiable with
 $$ \pddt u(t,x) = -\Delta u(t,x).$$
\end{proof}

\begin{theorem} 
   Let $G=(V,E)$ be a connected, locally finite graph. 
  Then the function $p: (0,\infty)\times V\times V\to\R_{\geq 0}$ is a fundamental solution for the 
  heat equation and does not depend on the choice of the exhaustion sequence $U_k$.   
\end{theorem}

\begin{proof} 
  The independence of $p$ from the choice of the exhaustion sequence follows from
  the maximum principle, more precisely from the domain monotonicity of $p_U$.\\
  To show that $p$ is a fundamental solution, we first remark that
  $p_k(t,x,y)\geq 0\, (x,y \in V),$ 
  $\sum_{y\in V}p_k(t,x,y) \leq 1\, (x\in V), $ and 
  $\pddt p_k(t,x,\cdot) = \Delta_{U_k} p_k(t,x,\cdot)$ for all $x$ in the interior of $U_k$. 
  By Proposition \ref{sequence of solutions}
  the sequence $p_k(t,x,y)$ converges for any $x\in V$ to a solution of the heat equation.\\
  Let now $u_0: V\to \R_{\geq 0}$ be a bounded, positive function (in the general case we split 
  the bounded function $u_0$ into its positive and negative part) and define
  $$ u_k(t,x) = \sum_{y\in V} p_k(t,x,y)u_0(y).$$
  We have
  \begin{eqnarray*}
   u_k(t,x)  &\leq &  \sup_{y\in V} u_0(y)\sum_{y\in V} p_k(t,x,y)\\
   			&\leq& \sup_{y\in V} u_0(y) 
  \end{eqnarray*}
  and hence, the sequence $u_k(t,\cdot)$ satisfies the assumptions (i) and (ii) in 
  Proposition \ref{sequence of solutions}. 
  As the sequence $u_k$ is non-decreasing its limit $u(t,x)=\lim_{k\to\infty} u_k(t,x)$
  is everywhere finite and satisfies the heat equation (cf. Proposition \ref{sequence of solutions}).\\
  Because of 
  $$ u(t,x) = \lim_{k\to\infty}\sum_{y\in V}p_k(t,x,y)u_0(y) = \sum_{y\in V}p(t,x,y)u_0(y)$$
  (note, that $p_k(t,x,y)$ is non-zero only for finitely many $y$)
  it remains to prove continuity at $t=0$, i.e. 
  $\lim_{t\to 0^+}u(t,x) = u_0(x).$
  
  To show this, we first prove that 
  $$ \lim_{t\to 0^+}\sum_{y\neq x}p(t,x,y) =0$$
  for any $x\in V$:
  If  $U\subset V$ is finite with $x\in \mathring{U}$ and $|\mathring{U}|=n+1$ we have
  $$   
  	1  \geq  \sum_{y\in V} p(t,x,y) \geq  p(t,x,x) \geq p_U(t,x,x) = 
		\sum_{j=0}^{n}\e^{-\lambda_j t} \Phi_j^2(x)  \to 
		\sum_{j=0}^n \Phi_j^2(x)\qquad (t\to 0^+).
  $$
  For any $x\in \mathring{U}$, the last sum equals one: if there was an $x\in \mathring{U}$
  such that $\sum_{j=0}^n \Phi_j^2(x)< 1$ we could conclude that
  $ \sum_{x\in \mathring{U}}\sum_{j=0}^n \Phi_j^2(x) < |\mathring{U}|=n+1$. But this would 
  contradict   $||\Phi_j|| = 1$.
 The claim now follows from
    $$
   1\geq \lim_{t\to 0^+}\sum_{y\in V} p(t,x,y) =
   	\lim_{t\to 0^+}\sum_{y \neq x}p(t,x,y) + \lim_{t\to 0^+}p(t,x,x).
   $$	 
  We therefore may conclude
  \begin{eqnarray*}
   \lim_{t\to 0^+} \big(u(t,x) - u_0(x)\big)  &=&
   		\lim_{t\to 0^+} \sum_{y\in V} p(t,x,y) \big(u_0(y) - u_0(x)\big)\\
		&=& \lim_{t\to 0^+} \sum_{y\neq x} p(t,x,y) \big(u_0(y) - u_0(x)\big).		 
  \end{eqnarray*}
  We obtain
  $$  \left| \sum_{y\neq x} p(t,x,y) \big(u_0(y) - u_0(x)\big)\right| \leq
  	2\sup u_0 \cdot \sum_{y\neq x} p(t,x,y) \longrightarrow 0\quad (t\to 0^+).$$
\end{proof}

It turns out that the heat kernel $p$ constructed above is the kernel of the 
heat semigroup $\e^{-t\bar{\Delta}}$:
\begin{theorem}\label{theorem heat kernel semigroup}
 For any $u_0\in C_c(V)$ we have
 $$ \e^{-t\bar{\Delta}}u_0(x) = \sum_{y\in V}p(t,x,y)u_0(y).$$
\end{theorem}
\noindent For the proof of this theorem, we will need the following two lemmas.
\begin{lemma}\label{lemma contraction}
The operator 
 \[
 P_t :\ell^2(V) \to \ell^2(V), \, P_tu(x) =  \sum_{y\in V}p(t,x,y)u(y)
 \]
 is a contraction for each $t\geq 0$.
\end{lemma}
\begin{proof}
Let $u_0\in C_c(V)$ and assume w.l.o.g. $u_0\geq 0$.
Choose $k\in\N$ large enough such that $\supp(u_0)\subset U_k$. Then we have (remember that
$\e^{-t\Delta_{U_k}}u_0 = \sum_{y\in U_k} p_k(t,\cdot,y)u_0(y)$)
 \begin{eqnarray*}
   || \sum_{y\in V} p_k(t,\cdot,y)u_0(y)||_{\ell^2(V)} &=& 
   					|| \sum_{y\in U_k} p_k(t,\cdot,y)u_0(y)||_{\ell^2(U_k)}\\
   		&=& || \e^{-t\Delta_{U_k}}u_0||_{\ell^2(U_k)} \\
		&\leq& ||u_0||_{\ell^2(U_k)}\\
		&=& ||u_0||_{\ell^2(V)}.
 \end{eqnarray*}
 This, together with Fatou's Lemma, yields
 \begin{eqnarray*}
  ||P_t u_0||^2_{\ell^2(V)}  &=& \sum_{x\in V}\Big| \sum_{y\in V} p(t,x,y)u_0(y)  \Big|^2\\
  		&\leq& \lim_{k\to\infty}  \sum_{x\in V}\Big| \sum_{y\in V} p_k(t,x,y)u_0(y)  \Big|^2\\
		&=& \lim_{k\to\infty}  ||\e^{-t\Delta_{U_k}}u_0||^2_{\ell^2(U_k)}\\
		&\leq& ||u_0||^2_{\ell^2(V)},
 \end{eqnarray*}
 in particular, $P_tu_0 \in \ell^2(V).$
\end{proof}
\begin{lemma}\label{lemma commuting}
For any $u_0\in C_c(V)$ and $t\geq 0$ we have
\[
 \Delta (P_t u_0) = P_t(\Delta u_0)
\] 
\end{lemma}
\begin{proof}
To see this, we remark that from Lemma \ref{heat kernel finite}
 it follows $\Delta_x p_k(t,x,y) = \Delta_y p_k(t,x,y)$ and that this formula also applies to the limit
 $p(t,x,y)$. By the self-adjointness of $\bar{\Delta}$ we obtain
 \begin{eqnarray*}
   \Delta (P_t u_0)(x)   &=& \Delta_x \sum_{y\in U_k}p(t,x,y)u_0(y) \\
   			      &=& \sum_{y\in U_k}\Big(\Delta_x p(t,x,y)\Big)u_0(y)\\
			      &=& \sum_{y\in U_k}\Big(\Delta_y p(t,x,y)\Big)u_0(y)\\
			      &=& \sum_{y\in U_k}p(t,x,y)\Delta_y u_0(y)\\
			      &=& P_t(\Delta u_0)(x) .
 \end{eqnarray*}

\end{proof}
\begin{proof}[Proof of Theorem \ref{theorem heat kernel semigroup}]
From Lemma \ref{lemma commuting} it follows that $\Delta (P_t u_0)\in \ell^2(V)$.
 This implies $P_tu_0 \in \D$ and therefore, the function
 $$ v(t,x) = P_t u_0(x) - \e^{-t\bar{\Delta}}u_0(x)$$
 is contained in $\D$, too. We are going to show that $v=0$:
 \begin{eqnarray*}
  \sum_{x\in V} v^2(t,x)  &=& \sum_{x\in V} \int_0^t \frac{\partial}{\partial \tau} v^2(\tau,x) d\tau\\
  				    &=& -2 \sum_{x\in V} \int_0^t  v(\tau,x)\bar{\Delta} v(\tau,x) d\tau\\
				    &=& -2\int_0^t \sum_{x\in V}v(\tau,x)\bar{\Delta} v(\tau,x) d\tau\\
				    &=& -2\int_0^t \langle v(\tau,\cdot), \bar{\Delta} v(\tau,\cdot)\rangle d\tau
				    		\leq 0.		
 \end{eqnarray*}
 as the Laplacian is positive and hence, it follows $v=0$. 
 The interchange of summation and integration in the calculation from above is justified by
 Tonelli's Theorem as (note that $P_t$ and $\e^{-t\bar{\Delta}}$ are contractions and 
 $\bar{\Delta}\e^{-t\bar{\Delta}} = \e^{-t\bar{\Delta}}\bar{\Delta}$)
 \begin{eqnarray*}
 \sum_{x\in V}|v(\tau,x)\bar{\Delta} v(\tau,x)|&\leq& ||v(\tau,\cdot)||\cdot ||\bar{\Delta} v(\tau,\cdot)||\\
 					&\leq & 2|| u_0|| \cdot  2||\bar{\Delta} u_0||
 \end{eqnarray*}
 and hence the ``iterated integrals'' are finite.
\end{proof}

\begin{corollary}
 The heat semigroup $\e^{-t\bar{\Delta}}$ is positive, i.e. $\e^{-t\bar{\Delta}}f\geq 0$ if
 $f\geq 0$. 
\end{corollary}

\subsection{Uniqueness of bounded solutions}

In this subsection we consider for a graph $G=(V,E)$  the Cauchy problem (CP)
\[ 
   \left\{
    \begin{array}{rcl}
     \pddt u + \Delta u &=& 0\\
      u(0,x) 			&=& u_0(x)
    \end{array}
    \right.
\]
on $[0,T)\times V$ with initial condition $u_0: V\to \R$.\\

A locally finite, connected graph $G=(V,E)$ admits a natural metric $d: V\times V\to \N$ that can be defined as follows. We define $d(x,x)=0$ for all $x\in V$. If $x\neq y$ there is a finite number of vertices 
$x=x_0\sim x_1\sim \ldots \sim x_k=y \in V$ that connect $x$ and $y$. Then $d(x,y)$
is the smallest number $k$ of such vertices.

\begin{theorem}\label{unique solution}
 Let $G=(V,E)$ denote a graph with the following property:
 there are $x_0\in V$ and $C\geq 0$ such that $\Delta d(\cdot, x_0) \geq -C$.
 Then a bounded solution $u$ of (CP) is uniquely determined by $u_0$.
\end{theorem}
\begin{proof}
 Let $M_1 = \sup\{ |u(t,x)| : t\in (0,T), x\in V\}, M_2 = \sup\{ |u_0(x)| : x\in V\}$ and consider for $R\in\N$ 
 the function
 $$ v(t,x) = u(t,x) - M_2 - \frac{M_1}{R}\Big(d(x,x_0) + Ct \Big).$$
 If we denote by 
 $$ B_R=B(x_0,R) = \{ x\in V : d(x,x_0)\leq R\}$$
 the ball with radius $R$ and center $x_0$, we always have 
 $\partial B_R \subset  \{ x\in V : d(x,x_0) =R\}$ and we may conclude
 $$ v(t,x) \leq 0$$
 if $(t,x) \in \Big(\{0\}\times B_R\Big)\cup \Big([0,T)\times \partial B_R\Big).$\\
 On $[0,T)\times \mathring{B}_R$ we have
 \begin{eqnarray*}
   \Big(\pddt + \Delta\Big)v(t,x)  &=& -\frac{M_1}{R}\Big(\Delta d(x,x_0) + C\Big) \leq 0.
 \end{eqnarray*}
 From the maximum principle it follows $v(t,x)\leq 0$ on $[0,T)\times B_R$ which is equivalent to
 $$ u(t,x) \leq M_2 + \frac{M_1}{R}\Big(d(x,x_0) + Ct \Big).$$
 Letting $R\to \infty$ we obtain $u(t,x) \leq M_2$ on $[0,T)\times V$. Repeating the arguments
 for $-u$ yields
 $$ |u(t,x)| \leq \sup_{x\in V} |u_0(x)|,$$
for $(t,x) \in [0,T)\times V$. The claim now follows by considering differences of bounded solutions
with same initial condition.
\end{proof}

\begin{corollary}
 Let $G=(V,E)$ denote a graph as in the theorem above. Then any bounded solution $u$ of
 (CP) with initial condition $u_0$ satisfies the inequality 
 $$ |u(t,x)| \leq  \sup_{x\in V} |u_0(x)|,$$
 for $(t,x) \in [0,T)\times V$.
\end{corollary}

The condition $\Delta d(\cdot,x_0) \geq -C$ from Theorem \ref{unique solution} is always satisfied if
the valence $x\mapsto m(x)$ is a bounded function on $V$. In this case, a proof of 
Theorem \ref{unique solution} can also be found in \cite{MR2246774}.\\
However, there are graphs with unbounded valence, such that this condition is fulfilled. 
In the following  example, we have $\Delta d(\cdot,x_0)=-2$.
\begin{figure}[htb] 
 \begin{center}
  \includegraphics{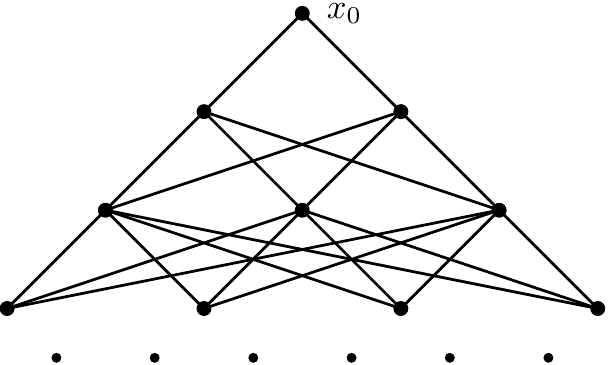}
  \end{center} 
   \caption{A Graph such that $\Delta d(\cdot,x_0) =-2.$} 
   \label{example} 
\end{figure}
We consider an infinite graph as in Figure \ref{example}. In the first row we have one
vertex $x_0$ which is connected to both vertices in the second row. In general, the $n$-th
row consists of $n$ vertices which are exactly connected to the $(n-1)$ vertices of the $(n-1)$-th row
and to the $(n+1)$ vertices of the $(n+1)$-th row. Then we obviously have
$$ \Delta d(x,x_0) = \sum_{y\sim x} \Big(d(x,x_0) - d(y,x_0)\Big) =-2.$$
Note, that the valence is unbounded.\\

In the smooth setting of complete Riemannian manifolds $M$ it was proved that bounded solutions
of the heat equation are unique if the Ricci curvature is bounded from below, see e.g. 
\cite{MR711862,MR505904}. On the other hand, a lower bound for the Ricci curvature of the
form $-(\dim(M)-1)\kappa^2$ implies the inequality
$$\Delta_M d(\cdot,x_0)\geq -(\dim(M)-1)\kappa\coth(\kappa d(\cdot,x_0)),$$ 
cf. \cite[Corollary I.1.2]{MR1333601}, and hence, our condition $\Delta d(\cdot,x_0) \geq -C$
from above can be interpreted as a weak curvature bound.

\begin{corollary}
 Let $G=(V,E)$ be a graph such that there are $x_0\in V$ and 
 $C\geq 0$ with $\Delta d(\cdot, x_0) \geq -C$.
 Then the following holds true:
 \begin{itemize}
  \item[\textup{(a)}] There exists a unique fundamental solution $p: (0,\infty)\times V\times V\to \R$
  		of the heat equation.
  \item[\textup{(b)}] $G$ is stochastically complete, i.e.
  		$$ \sum_{y\in V} p(t,x,y) =1$$
		for any $t>0$ and $x\in V$.
  \item[\textup{(c)}] For every $u_0\in \ell^1(V)$ and the corresponding bounded solution $u$ of (CP) we have
  		$$ \sum_{x\in V} u(t,x) = \sum_{x\in V} u_0(x)$$
		for any $t>0$.
 \end{itemize}
\end{corollary}
 \begin{proof}
  The claims in (a) and (b) follow immediately from the uniqueness of bounded solutions. To prove part  
  (c) we note that for the same reason we have
  $$ u(t,x) = \sum_{y\in V} p(t,x,y)u_0(y)$$
  and consequently
  $$ \sum_{x\in V} u(t,x) = \sum_{x\in V}\sum_{y\in V} p(t,x,y)u_0(y) 
  			    = \sum_{y\in V} u_0(y)\sum_{x\in V} p(t,x,y)
			    = \sum_{y\in V} u_0(y),$$
  where we used well known results on the rearrangement of absolutely convergent series and part
  (b). 	    		    
  \end{proof}
  
  It should be mentioned that Wojciechowski proved the equivalence of  stochastic completeness and
  the uniqueness of bounded solutions of the heat equation. Furthermore, he also showed
  that a locally finite graph is stochastically complete if and only if there is no 
  bounded, positive function  $f$ that satisfies the eigenvalue equation 
  $\Delta f = \lambda f$ for some $\lambda < 0$, cf. \cite[Theorem 3.2]{Wojciechowski:yq2008}.
  With this result at hand, Wojciechowski is able to provide examples for stochastically complete graphs
  which do not satisfy the condition  $\Delta d(\cdot, x_0) \geq -C$, 
  cf. \cite[Theorem 3.4]{Wojciechowski:yq2008} and the discussion following this Theorem.
  
\subsection*{Acknowledgements} 
I am deeply indebted to Daniel Lenz who carefully read the whole manuscript and made many very useful comments and suggestions.  I also want to thank Norbert Peyerimhoff for his interest in this work and for his helpful remarks. Furthermore, the unusually careful review of this paper is greatly 
acknowledged. \\


\bibliographystyle{amsplain}
\bibliography{dissertation,graphTheory}

\end{document}